\documentclass[letterpaper, 10 pt, conference]{ieeeconf}  

\IEEEoverridecommandlockouts                              

\usepackage{amsmath,amssymb,amsfonts}
\usepackage{algorithmic}
\usepackage{graphicx}
\usepackage{textcomp}
\usepackage{graphics} 
\usepackage{epsfig}
\usepackage{amsmath} 
\usepackage{amssymb}  
\usepackage{dsfont}
\usepackage{multirow}
\usepackage[dvipsnames]{xcolor}
\usepackage{flushend}

\usepackage[
    style=ieee,
    doi=false,
    isbn=false,
    url=true,
    eprint=false,
    backend=bibtex,
    maxcitenames=2,
    mincitenames=1,
    natbib=true
    ]{biblatex}
    
\bibliography{collection}

\usepackage{amsthm}
\newtheorem{theorem}{Theorem}[section]

\newtheorem{asmp}[theorem]{Assumption}

\newtheorem{df}[theorem]{Definition}
\newtheorem{exmp}[theorem]{Example}

\newcommand{\plim}[1]{\textup{p-lim}_{#1\to\infty}}
\newcommand{\T}{\intercal} 
\newcommand{\AL}{\underline{A}}
\newcommand{\AK}{\bar{A}}
\newcommand{\carsim}{CarSim\textsuperscript{\texttrademark}}
\newcommand{\spaceSave}{\textstyle}

\begin{document}

\title{Detecting Deception Attacks on Autonomous Vehicles\\ via Linear Time-Varying Dynamic Watermarking}

\author{Matthew~Porter$^{1}$,~Sidhartha~Dey$^{2}$,~Arnav~Joshi$^{1}$,~Pedro Hespanhol$^{3}$\\~Anil~Aswani$^{3}$,~Matthew~Johnson-Roberson$^{4}$,~and~Ram~Vasudevan$^{1}$
\thanks{*This work was supported by a grant from Ford Motor Company via the Ford-UM Alliance under award N022977 and by the UC Berkeley Center for Long-Term Cybersecurity.}%
\thanks{$^{1}$ M. Porter, A. Joshi, and R. Vasudevan are with the Department of Mechanical Engineering, University of Michigan, Ann Arbor, MI 48103, USA {\tt\small \{matthepo,arnavj,ramv\}@umich.edu}.}%
\thanks{S. Dey is with the Robotics Institute, University of Michigan, Ann Arbor, MI 48103, USA {\tt\small siddey@umich.edu}.}%
\thanks{P. Hespanhol and A. Aswani are with the Department of Industrial Engineering and Operations Research, University of California Berkeley, Berkeley, CA 94720, USA {\tt\small \{pedrohespanhol,aaswani\}@berkeley.edu}.}%
\thanks{M. Johnson-Roberson is with the Department of Naval Architecture, University of Michigan, Ann Arbor, MI 48103, USA {\tt\small mattjr@umich.edu}.}%
}

\maketitle

\begin{abstract}
Cyber-physical systems (CPS) such as autonomous vehicles rely on both on-board sensors and external communications to estimate their state. 
Unfortunately, these communications render the system vulnerable to cyber-attacks. 
While many attack detection methods have begun to address these concerns, they are limited to linear time-invariant (LTI) systems.
Though LTI system models provide accurate approximations for CPS such as autonomous vehicles at constant speed and turning radii, they are inaccurate for more complex motions such as lane changes, turns, and changes in velocity.
Since these more complex motions are more suitably described by linear time-varying (LTV) system models rather than LTI models, Dynamic Watermarking, which adds a private excitation to the input signal to validate measurements, has recently been extended to LTV systems. 
However, this extension does not allow for LTV systems that require several steps before the effect of a given control input can be seen in the measurement signal.
Additionally, there is no consideration for the time-varying effects of auto-correlation.
Furthermore, a proof of concept was only provided using simulations of a simplified model.

This paper relaxes the requirement for inputs to be visible in a single step and constructs an auto-correlation normalizing factor to remove the effects of auto-correlation. 
In addition, Dynamic Watermarking is applied to a high-fidelity vehicle model in \carsim and a 1/10 scale autonomous rover to further reinforce the proof of concept for realistic systems.
In each case, the vehicle follows a predefined path with time-varying velocity and turning radii.
A replay attack, which replays previously recorded measurements, is shown to be detectable using LTV Dynamic Watermarking in a quick and repeatable manner.
\end{abstract}
\section{Introduction}\label{sec:introduction}
For cyber-physical systems (CPS) to operate in a safe and efficient manner, their communications must remain secure. 
The difficulty of securing such systems has been illustrated in a variety of cases \cite{Langner2011,Abrams2008,Lee2016AnalysisGrid}. 
While most schemes  to detect whether such systems have been attacked have focused on linear time-invariant (LTI) systems, CPS such as autonomous vehicles (AV)s often require more complex models.
To overcome this, one particular scheme, Dynamic Watermarking, has recently been extended to linear time-varying (LTV) systems by Porter et al. \cite{Porter2019LTV}.
Unfortunately this extension does not allow for systems that require several steps before the effect of a given control input is seen in the measurement and does not take into account the time-varying effects of auto-correlation.
Furthermore, no real-world implementation of this method has been proposed.

\begin{figure}
    \centering
    \vspace{0.1in}
    \includegraphics[trim={0in 0in 0in 0in},clip,width=\columnwidth]{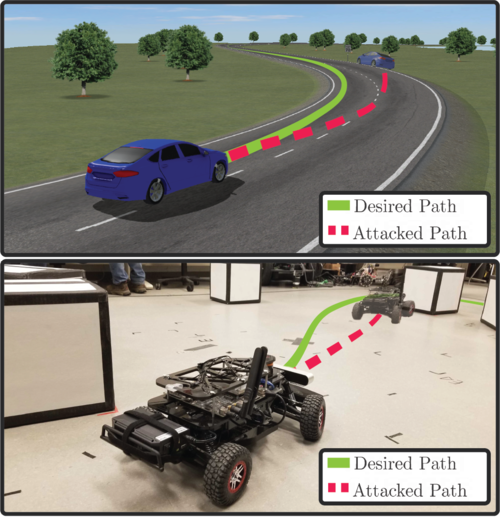}
    \caption{
    This paper describes the application and detection of deception attacks via LTV Dynamic Watermarking. 
    Proof of concept is shown using a high fidelity car model in \carsim (top) and a 1/10 scale autonomous rover (bottom). 
    }
    \label{fig:cover_pic}
\end{figure}

This paper relaxes the requirement for inputs to be visible in a single step and constructs an auto-correlation normalizing factor to remove the effects of auto-correlation resulting in a more consistent detection scheme.
In addition, this paper focuses on applying LTV Dynamic Watermarking to AVs performing complex motions both in high-fidelity simulation and in a real-world experiment as illustrated in Figure \ref{fig:cover_pic}. 

\subsection{Vulnerability of Autonomous Vehicles}
AV's have been touted as a way to increase safety by removing driver error.
However, like other CPS, AVs are vulnerable to cyber-attack, that give rise to additional safety concerns \cite{Gerla2015,Dominic2016,Amoozadeh2015}.
Some cyber-attacks can even compromise the control systems of the AV.
We divide these into two categories: \emph{direct attacks}, which seek to gain full control of the AV, and \emph{deception attacks}, which instead focus on attacking communications to alter the AV's perception of its surroundings.

\subsubsection{Direct Attacks}
AVs use on-board computers to handle sensor measurements and actuator inputs.
As a result, a hacker can take full control of the vehicle by compromising the security of these computers. 
Though modern cars have few avenues for accessing the on-board computer, they have been shown to be susceptible to hacking through wirelessly connecting to the infotainment system \cite{Greenberg2015} or through a wireless receiver connected to the diagnostic port \cite{Greenberg_2015}. 

\subsubsection{Deception Attacks} 
To detect obstacles and localize themselves, AVs use some combination of cameras, light detection and ranging (LiDAR), and GPS each of which has their own vulnerabilities.
While cameras are susceptible to glare \cite{petit2015remote}, further vulnerabilities lie in deceiving the object classifier that is run on the resulting images \cite{Eykholt2018}.
Moreover, fabricated LiDAR returns can be injected via lasers in an efficient enough manner to fool object detection algorithms \cite{Cao2019}.
Lastly, GPS measurements are susceptible to spoofing attacks with commercially available hardware \cite{Humphreys2008}.

\subsection{Attack Detection Algorithms}
This paper focuses on a particular detection method, Dynamic Watermarking, which adds a watermark signal to the control input to secure the measurement signals \cite{Mo2009,Satchidanandan2017,Hespanhol2017,Hespanhol2018,hespanhol2019sensor,Porter2019LTV}.
Dynamic Watermarking is an \emph{active method} meaning it makes alterations to the system which allow it to detect more complex attacks \cite{Porter2019,Weerakkody2016info,Weerakkody2017exposing}.
While, the specific methods of Dynamic Watermarking vary, in one popular approach the watermark is a multivariate Gaussian sequence that is generated by the controller.
These works present their detection schemes using two forms of tests: \emph{asymptotic tests}, which provide guarantees of detection in infinite time, and \emph{statistical tests}, which use the motivation of the asymptotic tests to form an implementable real-time attack detection scheme.
Both forms use the \emph{measurement residual}, defined as the difference between the measurement and the expected measurement, by considering its   covariance and correlation with the watermark.

Recently, Dynamic Watermarking was extended to LTV systems by Porter et al. \cite{Porter2019LTV}.
However, this extension requires the effect of the control input to be visible in the measurement signal  immediately. 
As a result, systems that do not fit this requirement such as those with distributed sensing and control are unable to use this extension.
Additionally there is no consideration for the effect of auto-correlation of the measurement residual.
While auto-correlation remains consistent for LTI Dynamic Watermarking, there is no guarantee of consistency for the LTV case.
Furthermore, this extension only provides proof of concept using a kinematic vehicle model with added Gaussian noise.

\subsection{Contributions}
The contributions of this paper are four-fold. 
First, a relaxed version of the requirement for inputs to be visible in a single step is provided. 
Using the less constrictive requirement, we show that the resulting asymptotic tests maintain the same detection guarantees.
Second, an auto-correlation normalizing factor is derived for use in the statistical tests to ensure that the resulting test metric is unaffected by the time-varying auto-correlation of the measurement residual.
Third, LTV Dynamic Watermarking is applied to a high fidelity car model in \carsim.
Fourth, LTV Dynamic Watermarking is applied to a 1/10 scale autonomous rover.
In each case the vehicle performs a path following task for a predefined path with time-varying velocity and turning radii.
A replay attack, which replays previously recorded measurements, is then shown to be detectable using LTV Dynamic Watermarking.
The simulation and the experiment are repeated several times to show the consistency of detection.

The remainder of this paper is organized as follows.
Section \ref{sec:Notation} covers the notation used throughout the paper.
Section \ref{sec:LTV_System} defines the LTV system and attack model.
Sections \ref{sec:asymp} provides an overview of the asymptotic tests developed by Porter et al. \cite{Porter2019LTV} and the proposed asymptotic tests.
Sections \ref{sec:stat} provides an overview of the statistical test developed by Porter et al. \cite{Porter2019LTV} and the proposed statistical test.
Sections \ref{sec:sim} and \ref{sec:experimental} validate these methods through simulated experiments in \carsim and physical experiments on a 1/10 scale autonomous rover, respectively.

\section{Notation}\label{sec:Notation}
The 2-norm of a vector $x$ is denoted $\|x\|$.
Similarly, the 2-norm of a matrix $X$ is denoted $\|X\|$.
The trace of a matrix $X$ is denoted tr$(X)$.
Zeros matrices of dimension $i\times j$ are denoted $0_{i\times j}$ and in the case that $i=j$ the notation is simplified to $0_i$.
Identity matrices of dimension $i$ are denoted $I_i$.

The Wishart distribution with scale matrix $\Sigma$ and $i$ degrees of freedom is denoted $\mathcal{W}(\Sigma,i)$ \cite[Section 7.2]{anderson2003}.
The multivariate Gaussian distribution with mean $\mu$  and covariance $\Sigma$ is denoted $\mathcal{N}(\mu,\Sigma)$.
The matrix Gaussian distribution with mean $\mathcal{M}$, and parameters $\Sigma$ and $\Omega$ is denoted $\mathcal{N}(\mathcal{M},\Sigma,\Omega)$.
The expectation of a random variable $a$ is denoted $\mathds{E}[a]$.
Given a sequence of random variables $\{a_i\}_{i=1}^\infty$, convergence in probability is denoted $\text{p-lim}_{i\to\infty}a_i$ 
\cite[Definition 7.2.1]{grimmett2001probability}.

\section{LTV System and Attack Model}\label{sec:LTV_System}
This section provides an overview of the LTV system and attack model used by Porter et al. \cite{Porter2019LTV}.
Consider an LTV system with state $x_n$, measurement $y_n$, process noise $w_n$, measurement noise $z_n$, watermark $e_n$, additive attack $v_n$, and stabilizing feedback that uses the observed state $\hat{x}_n$
\begin{align}
    x_{n+1}&=A_nx_n+B_nK_n\hat{x}_n+B_ne_n+w_n\label{eq:state_update}\\
    y_n&=C_nx_n+z_n+v_n\label{eq:output_equation}
\end{align}
where $x_n,\hat{x}_n,w_n\in \mathbb{R}^p$, $e_n\in\mathbb{R}^q$, $y_n,z_n,v_n\in\mathbb{R}^r$, and $x_0=0_{p\times 1}$.
Note, the initial condition of zero is not necessary since the effects of a non-zero initial condition would asymptotically decay under the assumption of stability.
However, this assumption eases notation.
The process noise $w_n$, measurement noise $z_n$, and watermark $e_n$ are mutually independent and take the form $w_n\sim\mathcal{N}(0_{p\times 1},\Sigma_{w,n})$, $z_n\sim\mathcal{N}(0_{r\times 1},\Sigma_{z,n})$, and $e_n\sim\mathcal{N}(0_{q\times 1},\Sigma_e)$. 
While the process and measurement noise are unknown to the controller, the watermark signal is generated by the controller and is known.
For simplicity, define $\AK_n=(A_n+B_nK_n)$ and $\AK_{(n,m)}=\AK_n\cdots\AK_m$ for $n\geq m$ and $\AK_{(n,n+1)}=I_p$.
Then make the following assumption.
\begin{asmp}\label{asmp:bound1}
The covariances $\Sigma_e$, $\Sigma_{w,n}$, and $\Sigma_{z,n}$, of the random variables used in \eqref{eq:state_update}-\eqref{eq:output_equation}, are full rank.
Furthermore, there exists positive constants $\eta_w,\eta_z,\eta_{\AK},\eta_B,\eta_C\in\mathbb{R}$ such that $\|\Sigma_{w,n}\|<\eta_w$, $\|\Sigma_{z,n}\|<\eta_z$, $\|\AK_n\|<\eta_{\AK}<1$, $\|B_n\|<\eta_B$, and $\|C_n\|<\eta_C$, for all $n\in\mathbb{N}$.
\end{asmp}
\noindent The assumption of bounded full rank covariances for the process and measurement noise are satisfied for most systems by modeling error and sensor noise.
Furthermore, the input and output matrices are often constrained to be finite by sensor and actuator limits.
Since the watermark covariance is user defined, the assumption of full rank can be satisfied by choosing a positive definite matrix.
Additionally, since the controller is user defined, the constraint on $\AK_n$ can often be satisfied by proper choice of $K_n$.
Then make the following assumption.
\begin{asmp}\label{asmp:1step}
\begin{align}
    \textstyle\lim_{i\to\infty}~\spaceSave\frac{1}{i}\sum_{n=0}^{i-1} C_nB_{n-1}\neq 0_{r\times q}.\label{eq:consistently_observervable_condition}
\end{align}
\end{asmp}
\noindent Here, \eqref{eq:consistently_observervable_condition} guarantees a persistent correlation between the measurement signal $y_n$ and the watermark signal $e_{n-1}$, which has been delayed by a single time step.
This ensures that the watermark has a persistent measurable effect on the measurement signal, which can then be used for validation purposes.
This assumption is later replaced in the proposed tests.

The observer and the corresponding observer error, defined as $\delta_n=\hat{x}_n-x_n$, satisfy 
\begin{align}
    \hat{x}_{n+1}&=(\AK_n+L_nC_n)\hat{x}_n+B_ne_n-L_ny_n\label{eq:obsver_update}\\
    \delta_{n+1}&=(A_n+L_nC_n)\delta_n-w_n-L_n(z_n+v_n),\label{eq:observer_error_update_full}
\end{align}
where $\hat{x}_0=\delta_0=0_{p\times 1}$. 
For simplicity, define $\AL_n=(A_n+L_nC_n)$ and $\AL_{(n,m)}=\AL_n\cdots\AL_m$ for $n\geq m$ and $\AL_{(n,n+1)}=I_p$.

Next, consider the expected value $\Sigma_{\delta,n}=\mathds{E}[\delta_n\delta_n^\T~|~v_n=0_{r\times 1},~\forall n]$, which can be written as
\begin{align}
    \textstyle\Sigma_{\delta,n}&=\spaceSave\sum _{i=0}^n \AL_{(n-1,n-i+1)}(\Sigma_{w,n-i}+\nonumber\\
    &\qquad+L_{n-i}\Sigma_{z,n-i}L_{n-i}^\T)\AL_{(n-1,n-i+1)}^\T.\label{eq:observer_error_covariance_normal}
\end{align}
The \emph{matrix normalization factor} is then defined as
\begin{align}
    V_n=(C_n\Sigma_{\delta,n}C_n^\T+\Sigma_{z,n})^{-1/2}\label{eq:residual_normalizer},
\end{align}
which exists since $\Sigma_{z,n}$ is full rank. 
For the LTV system, the matrix normalization factor can be thought of as a time-varying normalization for the covariance of the measurement residual.
Next, make the following assumption about the observer.
\begin{asmp}\label{asmp:obs_bounds}
    There exists positive constants $\eta_{\AL},$ $\eta_L,$ $\eta_\delta,$
    $\eta_V\in\mathbb{R}$ such that $\|\AL_n\|<\eta_{\AL}<1$, $\|L_n\|<\eta_L$, $\|\Sigma_{\delta,n}\|<\eta_\delta$, and $\|V_n\|<\eta_V$, for all $n\in\mathbb{N}$.
\end{asmp}
\noindent Since the observer is user defined, the assumptions on $\AL_n$ and $L_n$ can often be satisfied for proper choice of $L_n$. 
The bound on $\Sigma_{\delta,n}$ is satisfied since applying previous assumptions to \eqref{eq:observer_error_covariance_normal} results in an increasing geometric series with finite bound. 
The assumption on $V_n$ is satisfied by lower bounding the eigenvalues of the measurement noise which is often satisfied.

Next, consider an attack $v_n$ that satisfies
\begin{align}
    v_n&=\alpha(C_nx_n+z_n)+C_n\xi_n+\zeta_n\label{eq:attack_definition}\\
    \xi_{n+1}&=\AK_n\xi_n+\omega_n\label{eq:false_state_update},
\end{align}
where $\alpha\in\mathbb{R}$ is called the \emph{attack scaling factor}, the \emph{false state} $\xi_n\in\mathbb{R}^p$ has process noise $\omega_n\in\mathbb{R}^p$ and measurement noise $\zeta_n\in\mathbb{R}^r$ that take the form $\omega_n\sim\mathcal{N}(0_{p\times 1},\Sigma_{\omega,n})$ and  $\zeta_n\sim\mathcal{N}(0_{r\times 1},\Sigma_{\zeta,n})$ and are mutually independent with each other and with $w_n$ and $z_n$. 
When $\Sigma_{\omega,n}$ and $\Sigma_{\zeta,n}$ are selected properly and the attack scaling parameter is $-1$, this model describes a replay attack.
While an attacker could choose to allow the noise to have unbounded covariance, the resulting attack would be trivial to detect.
Therefore, make the following assumption about the attack model.
\begin{asmp}
When there is an attack, $v_n$ follows the dynamics \eqref{eq:attack_definition}-\eqref{eq:false_state_update} with the attack scaling factor remaining constant.
Furthermore, there exists positive constants $\eta_\omega,\eta_\eta\in\mathbb{R}$ such that $\|\Sigma_{\omega,n}\|<\eta_\omega,~\|\Sigma_{\zeta,n}\|<\eta_\zeta$, for all $n\in\mathbb{N}$.
\end{asmp}

Though in a real-world attack the attacker could likely start and stop the attack as desired, attacks that are not consistently present are impossible to detect every time due to the noise in the system.
As a result, some notion of the persistence must be defined to make asymptotic guarantees of detection.
LTV Dynamic Watermarking uses the following definition to describe the persistence of an attack.
\begin{df}\label{def:LTV_power}
The \underline{asymptotic attack power} is defined as
\begin{align}
    \spaceSave\plim{i}~\frac{1}{i}\sum_{n=0}^{i-1}v_n^\T v_n.\label{eq:attack_power}
\end{align}
\end{df}
\noindent Here, an asymptotic attack power greater than 0 is considered to be a persistent attack.

\section{Asymptotic Tests}\label{sec:asymp}
This section describes the asymptotic tests developed by Porter et al. \cite{Porter2019LTV} and the proposed asymptotic tests which include a relaxed form of Assumption \ref{asmp:1step}.
In each case, the tests are proven to detect the generalized replay attack described in \eqref{eq:attack_definition}-\eqref{eq:false_state_update} in infinite time.

\subsection{Previous Asymptotic Tests}\label{sec:asymp_old}
The asymptotic guarantee of detection for the tests developed by Porter et al. \cite{Porter2019LTV} are defined in the following theorem.
\begin{theorem}\cite[Theorem III.6]{Porter2019LTV}\label{thm:LTV_Asymptotic_Main_Result}
Consider an attacked LTV system satisfying the dynamics in \eqref{eq:state_update}-\eqref{eq:output_equation} and \eqref{eq:obsver_update}-\eqref{eq:observer_error_update_full}. 
Let $V_n$ be as defined in \eqref{eq:residual_normalizer}.
If $v_n=0_{r\times 1}~\forall n$, then 
\begin{align}
    \plim{i}~&\spaceSave\frac{1}{i}\sum_{n=0}^{i-1} V_n(C_n\hat{x}_n-y_n)e_{n-1}^\T=0_{r\times q} \tag{C1} \label{eq:ltv_watermark_correlation_test}\\
    \intertext{and}
    \plim{i}~&\spaceSave\frac{1}{i}\sum_{n=0}^{i-1} V_n(C_n\hat{x}_n-y_n)(C_n\hat{x}_n-y_n)^\T V_n^\T=I_r \tag{C2} \label{eq:ltv_covariance_test}.
\end{align}
Furthermore, if the attack follows the dynamics in \eqref{eq:attack_definition}-\eqref{eq:false_state_update} and has non-zero asymptotic power as defined in \eqref{eq:attack_power}, then \eqref{eq:ltv_watermark_correlation_test} and \eqref{eq:ltv_covariance_test} cannot both be satisfied.
\end{theorem}
\noindent Note, the LHS of \eqref{eq:ltv_watermark_correlation_test} and \eqref{eq:ltv_covariance_test} can be used to guarantee detection of generalized replay attacks with non-zero asymptotic power in infinite time.
Since the statistical tests consider only a finite number of steps at a time, the sample averages of the measurement residuals covariance and correlation with the watermark are more likely to be closer to the RHS of \eqref{eq:ltv_watermark_correlation_test} and \eqref{eq:ltv_covariance_test} when no attack is present. 
As a result, the test becomes more sensitive.

\subsection{Proposed Asymptotic Tests}\label{sec:asymp_mod}
To relax the requirement made in Assumption \ref{asmp:1step}, we instead consider the following assumption. 
\begin{asmp}\label{asmp:delay}
There exists $\kappa\in\mathbb{N}$ such that 
\begin{align}
    \spaceSave\lim_{i\to\infty}~\frac{1}{i}\sum_{n=0}^{i-1} C_n\AK_{(n-1,n-\kappa+1)}B_{n-\kappa}\neq 0_{r\times q}.\label{eq:delay_observervable_condition}
\end{align}
\end{asmp}
\noindent Here \eqref{eq:delay_observervable_condition} guarantees a persistent correlation between the measurement signal $y_n$ and the watermark signal $e_{n-\kappa}$ which has been delayed by $\kappa$ time steps. 
This allows for systems that require more than a single step for the watermark to have a persistent measurable effect on the measurement signal.

Next we replace Theorem \ref{thm:LTV_Asymptotic_Main_Result} with the following theorem.
\begin{theorem}\label{thm:new_LTV_Asymptotic_Main_Result}
Consider an attacked LTV system satisfying the dynamics in \eqref{eq:state_update}-\eqref{eq:output_equation} and \eqref{eq:obsver_update}-\eqref{eq:observer_error_update_full}. 
Let $V_n$ be as defined in \eqref{eq:residual_normalizer} and $\kappa$ be the smallest value for which \eqref{eq:delay_observervable_condition} holds.
If $v_n=0_{r\times 1}~\forall n$, then 
\begin{align}
    \plim{i}~&\spaceSave\frac{1}{i}\sum_{n=0}^{i-1} V_n(C_n\hat{x}_n-y_n)e_{n-\kappa}^\T=0_{r\times q} \tag{$\mathcal{C}$1} \label{eq:new_ltv_watermark_correlation_test}\\
    \intertext{and}
    \plim{i}~&\spaceSave\frac{1}{i}\sum_{n=0}^{i-1} V_n(C_n\hat{x}_n-y_n)(C_n\hat{x}_n-y_n)^\T V_n^\T=I_r \tag{$\mathcal{C}$2} \label{eq:new_ltv_covariance_test}.
\end{align}
Furthermore, if the attack follows the dynamics in \eqref{eq:attack_definition}-\eqref{eq:false_state_update} and has non-zero asymptotic power as defined in \eqref{eq:attack_power}, then \eqref{eq:new_ltv_watermark_correlation_test} and \eqref{eq:new_ltv_covariance_test} cannot both be satisfied.
\end{theorem}
\noindent Note, the delay of $e_{n-\kappa}$ in \eqref{eq:new_ltv_watermark_correlation_test} has been changed but the guarantees of the theorem remain the same.
The proof of Theorem \ref{thm:new_LTV_Asymptotic_Main_Result} follows the same format as the proof of Theorem \ref{thm:LTV_Asymptotic_Main_Result} but uses the following theorem in place of \cite[Theorem III.7]{Porter2019LTV}.

\begin{theorem}\label{thm:alpha}
Consider an attacked LTV system satisfying the dynamics in \eqref{eq:state_update}-\eqref{eq:output_equation} and \eqref{eq:obsver_update}-\eqref{eq:observer_error_update_full} and an attack model satisfying \eqref{eq:attack_definition}-\eqref{eq:false_state_update}.
Let $V_n$ be as defined in \eqref{eq:residual_normalizer}, and $\kappa$ being the smallest value for which \eqref{eq:delay_observervable_condition} holds. \eqref{eq:new_ltv_watermark_correlation_test} holds if and only if the attack scaling factor $\alpha$ is equal to 0.
\end{theorem}
\noindent The proof for Theorem \ref{thm:alpha} is found in the appendix.

\begin{figure*}[t]
    \centering
    \includegraphics[trim={0.45in 0 2.65in 0},clip,height=2.05in]{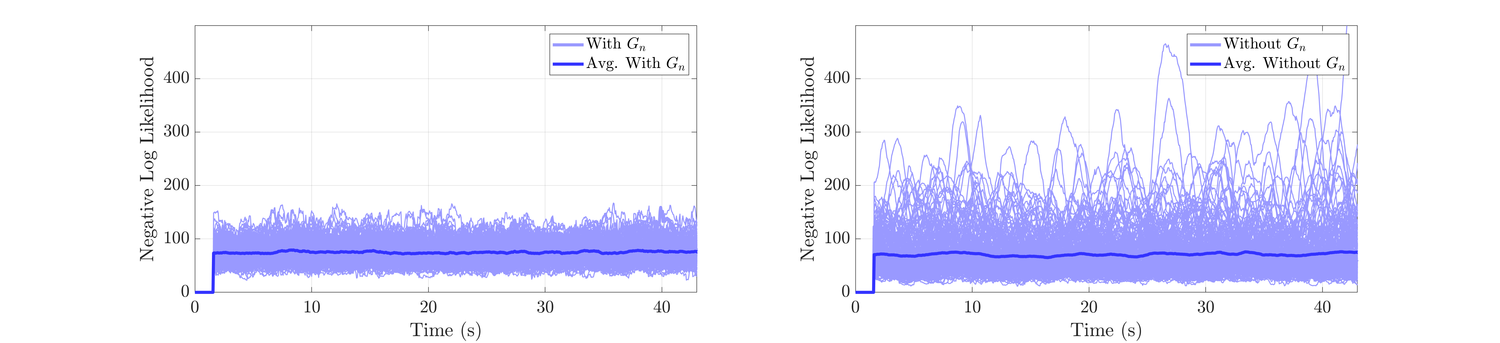}
    \includegraphics[trim={2.6in 0 0.45in 0},clip,height=2.05in]{figs/G_example2_lq.png}
    \caption{ The LTV system in Example \ref{exmp:auto} is simulated 200 times and the negative log likelihood is generated with the auto-correlation normalizing factor, $G_n$, (left) and without the auto-correlation normalizing factor (right).
    }
    \label{fig:G_compare}
\end{figure*}
 
\section{Statistical Tests}\label{sec:stat}
This section describes the statistical test developed by Porter et al. \cite{Porter2019LTV} and the proposed statistical test which uses Assumption \ref{asmp:delay} in place of Assumption \ref{asmp:1step} and the auto-correlation normalizing factor.
The addition of this normalizing factor is then shown to provide a more consistent test metric for an example LTV system.

\subsection{Previous Statistical Test}\label{sec:stat_old}
While Section \ref{sec:asymp} provides necessary background for LTV Dynamic Watermarking, infinite limits are not well suited for real time attack detection. 
This section derives a statistical test using a sliding window approach.
Let 
\begin{align}
    \psi_n&=\begin{bmatrix}V_n(C_n\hat{x}_n-y_n)\\ e_{n-1}\end{bmatrix}\\\intertext{and}
    Q_n&=[\psi_{n-\ell}~\hdots~\psi_{n}][\psi_{n-\ell}~\hdots~\psi_{n}]^\T.
\end{align}
where $\ell+1$ is the window size, $\ell\in\mathbb{N}$, and $\ell\geq q+r-1$.
Note, $\psi_n$ is asymptotically uncorrelated and identically distributed such that $\psi_n\sim\mathcal{N}(0_{q+r\times1},S)$, for $n=1,2,3,\cdots$ where
\begin{align}\label{eq:stat_cov}
    S=\begin{bmatrix}I_r & 0_{r\times q}\\0_{q\times r} & \Sigma_e\end{bmatrix}.
\end{align}
Therefore, under the assumption of no attack, the distribution of $Q_n$ approaches a Wishart distribution with $\ell+1$ degrees of freedom and scale matrix $S$ as $\ell$ goes to infinity.
Furthermore, for a generalized replay attack with non-zero asymptotic power, Theorem \ref{thm:LTV_Asymptotic_Main_Result} proves that the scale matrix for $Q_n$ is no longer $S$ since either \eqref{eq:ltv_watermark_correlation_test} or \eqref{eq:ltv_covariance_test} is not satisfied. 
The Wishart distribution can then be used to define a statistical test using the negative log likelihood of the scale matrix $S$ given the sampled matrix $Q_n$:
\begin{align}\label{eq:nll}
\mathcal{L}(Q_n)=(q+r-\ell)\log(|Q_n|)+tr(S^{-1}Q_n).
\end{align}
For a user defined threshold, a negative log likelihood greater than the threshold raises an alarm.

In theory, if the process and measurement noise covariances $\Sigma_{w,n}$ and $\Sigma_{z,n}$ are known, $V_n$ can be calculated using \eqref{eq:observer_error_covariance_normal}-\eqref{eq:residual_normalizer}.
In practice, these covariances are difficult to estimate which can lead to error in the estimate of $V_n$. 
To reduce this error, $V_n$ can be directly estimated using the ensemble average  
\begin{align}\label{eq:vn_approx}
V_n\approx\left(\spaceSave\frac{1}{i}\sum_{j=1}^{i}(C_n\hat{x}_n^{(j)}-y_n^{(j)})(C_n\hat{x}_n^{(j)}-y_n^{(j)})^\T\right)^{-1/2}
\end{align}
where the superscript $(j)$ is the index of the realization and the approximation improves as $i$ becomes larger.
This approximation is appropriate since by the weak law of large numbers we have that when no attack is present
\begin{align}
    \plim{i}~\spaceSave\frac{1}{i}\sum_{j=1}^{i}&(C_n\hat{x}_n^{(j)}-y_n^{(j)})(C_n\hat{x}_n^{(j)}-y_n^{(j)})^\T=\nonumber\\
    &=C_n\Sigma_{\delta,n}C_n^\T+\Sigma_{z,n}
\end{align}
and $V_n$ is defined as in \eqref{eq:residual_normalizer}.

\begin{figure*}
    \centering
    \includegraphics[trim={0.55in 0 2.7in 0},clip,height=2.1in]{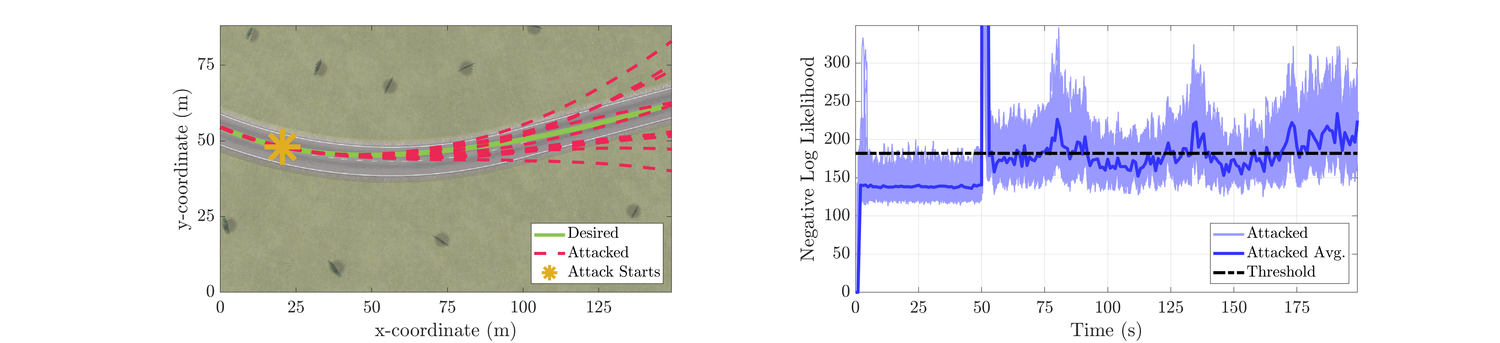}
    \includegraphics[trim={2.6in 0 0.45in 0},clip,height=2.1in]{figs/carsim_fig_lq.png}
    \caption{
    The simulated high fidelity car is attacked with a replay attack after 50 s of operation.
    The desired trajectory and 10 attacked realizations are plotted for the region that the attack is initiated (left). 
    Negative log likelihood for all 200 attacked realizations with average value are plotted (right).
    }
    \label{fig:carsim}
\end{figure*}

\subsection{Modified Statistical Test}\label{sec:stat_mod}
By adding an auto-correlation normalizing factor, denoted $G_n$, and the delay of $\kappa$ to the watermark, the proposed statistical test is as follows.
Let 
\begin{align}
    \psi_n=\begin{bmatrix}
    V_n(C_n\hat{x}_n-y_n)\\ e_{n-\kappa}
    \end{bmatrix}~\text{and}~
    P_n=\begin{bmatrix}
    \psi_{n-\ell}~\hdots~\psi_{n}
    \end{bmatrix},
\end{align}
where $\ell+1$ is the window size, $\ell\in\mathbb{N}$, and $\ell\geq q+r-1$.
Then $P_n$ is distributed according to
\begin{align}
P_n&\sim\mathcal{N}(0_{(\ell+1)\times(q+r)},S,G_n),
\intertext{where $S$ is as defined in \eqref{eq:stat_cov} and }
G_n&=\spaceSave\frac{\mathds{E}[P_n^\T P_n]}{\text{tr}(S)}.
\label{eq:scale_autocorrelation_matrices}
\end{align}
We can calculate the elements of $G_n$ as 
\begin{align}
&\spaceSave\frac{\mathds{E}[\psi_{j+i}^\T\psi_{j}]}{\text{tr}(S)}=\frac{\mathds{E}[\psi_{j}^\T\psi_{j+i}]}{\text{tr}(S)}=\nonumber\\
&\spaceSave=\frac{\text{tr}(C_j^\T V_j^\T V_{j+i}C_{j+i}\AL_{j+i-1}\hdots\AL_{j}\Sigma_{\delta,j})+\text{tr}(\Sigma_e)}{\text{tr}(S)}+\nonumber\\
&\hspace{0.5in}+\spaceSave\frac{\text{tr}(V_j^\T V_{j+i}C_{j+i}\AL_{j+i-1}\hdots\AL_{j+1}L_j\Sigma_{z,j}
)}{\text{tr}(S)}\label{eq:autocorrelation_matrix_elements}.
\end{align}
Finally, by \citet[Theorem 2.4.1]{Kollo2005} we have that
\begin{align}
    Q_n=P_n G_n^{-1} P_n^\T\sim\mathcal{W}_{q+r}(S,\ell+1).
\end{align}
Note, due to the addition of the auto-correlation normalizing factor, $Q_n$ is distributed according to a Wishart distribution for all $\ell\geq q+r-1$ instead of only approaching a Wishart distribution.
Furthermore, for a generalized replay attack of non-zero asymptotic power, Theorem \ref{thm:new_LTV_Asymptotic_Main_Result} proves that the scale matrix for $Q_n$ is no longer $S$ since either \eqref{eq:new_ltv_watermark_correlation_test} or \eqref{eq:new_ltv_covariance_test} is not satisfied.
Using this new definition of $Q_n$ the statistical test follows that of \eqref{eq:nll}.
Note, Allowing $G_n$ to be $I_{(\ell+1)}$ for all $n\in\mathbb{N}$ and assuming $\kappa=1$ results in the statistical test provided by Porter et al. \cite{Porter2019LTV}.

To avoid compounding error from estimated process and measurement noise, $G_n$ can be directly estimated using the ensemble average
\begin{align}
G_n\approx\spaceSave\frac{1}{i}\sum_{j=1}^i \frac{\left(P_n^{(j)}\right)^\T P_n^{(j)}}{\text{tr}(S)},\label{eq:gn_approx}
\end{align}
where the superscript $(j)$ is the index of the realization and the approximation improve as $i$ becomes larger. 
This approximation is appropriate since when no attack is present
\begin{align}
     \plim{i}~\spaceSave\frac{1}{i}\sum_{j=1}^{i}\frac{\psi_{n}^{(j)\T}\psi_{n+k}^{(j)}}{\text{tr}(S)}=\frac{\mathds{E}[\psi_{n}^{\T}\psi_{n+k}]}{\text{tr}(S)}.
\end{align}

To illustrate the effect of the auto-correlation normalizing factor consider the following example.
\begin{exmp}\label{exmp:auto}
Consider an LTV system satisfying the dynamics in \eqref{eq:state_update}-\eqref{eq:output_equation} where $v_n=0$ for all $n$, $w_n\sim\mathcal{N}(0_{3\times1},1\times10^{-3}I_3)$, $z_n\sim\mathcal{N}(0_{2\times1},1\times10^{-3}I_2)$, $e_n\sim\mathcal{N}(0,1\times10^{-3})$,
\begin{align}
    A_n&=\begin{bmatrix}
    1 & 1+\frac{1}{2}\sin(\frac{n}{100})&0\\
    0&1&0.1\\
    0&0&1
    \end{bmatrix},\\
    B_n&=\begin{bmatrix}0&0&1\end{bmatrix}^\T,\\
    C_n&=\begin{bmatrix}1 & 0 & 0\\0 & 1 & 0\end{bmatrix},
\end{align}
and
\begin{align}
    K_n&=\begin{bmatrix}-4\times10^{-4} & -3.65\times10^{-2} & -1.05\times10^{-1}\end{bmatrix}.
\end{align}
Furthermore consider an observer satisfying \eqref{eq:obsver_update} where
\begin{align}
    L_n=\begin{bmatrix}-7\times10^{-2} & -1\\-2.2\times10^{-3} & -1.4\times10^{-1}\\-1.6\times10^{-3} & -4.5\times10^{-2}\end{bmatrix}.
\end{align}
Note, for this system $\kappa$ is 2. 
The test metric was generated for 200 simulated realizations both with and without the auto-correlation normalizing factor $G_n$ for a window size of 20 ($\ell=19$).
\end{exmp}
As illustrated in Figure \ref{fig:G_compare}, the addition of the auto-correlation normalizing factor has little effect on the average of the negative log likelihood. 
However, this normalizing factor does improve the consistency by removing anomalies in many of the realizations caused by auto-correlation.

\begin{figure*}[t]
    \centering
    \includegraphics[trim={0.7in 0 2.8in 0},clip,height=2.15in]{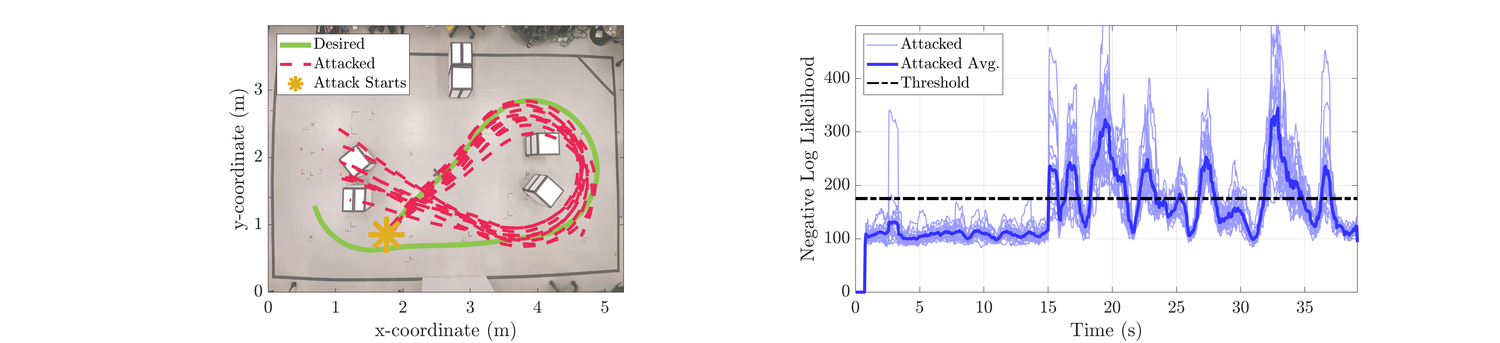}
    \includegraphics[trim={2.5in 0 0.4in 0},clip,height=2.15in]{figs/rover_fig_lq.png}
    \caption{
    The 1/10 scale autonomous rover is attacked with a replay attack after 15 s of operation. 
    The desired trajectory and 10 attacked realizations are plotted for the region that the attack is initiated (left). 
    Negative log likelihood for all 20 attacked realizations with average value are plotted (right).
    }
    \label{fig:rover}
\end{figure*}

\section{Simulated Results}\label{sec:sim}
This section illustrates the effectiveness of LTV Dynamic Watermarking using a high fidelity vehicle model in \carsim. 
For the simulation, the vehicle completes a 1,137 m long trajectory traveling at speeds up to 7 m/s in approximately 200 s.
This is accomplished using a linear quadratic regulator (LQR) and a linearization of the car model \cite{Shuai2014}.
The simulated measurement signal at step $n$ includes the ground plane coordinates $(x_{c,n},y_{c,n})$ in meters, heading $\psi_n$ in radians, longitudinal velocity $v_{1,n}$ in meters per second, lateral velocity $v_{2,n}$ in meters per second, yaw rate $\dot{\psi}_n$ in radians per second, and steering wheel angle $\delta_n$ in radians.
Since the feedback from the simulation does not include noise, Gaussian measurement noise was added to the measurement such that when no attack is present
\begin{align}
    y_n=\begin{bmatrix}
    x_{c,n} & y_{c,n} & \psi_n & v_{1,n} & v_{2,n} & \dot{\psi}_n & \delta_n
    \end{bmatrix}^\T+z_n,
\end{align}
where 
\begin{align}
    z_n\sim\mathcal{N}\left(0_{7\times1},
    1\times10^{-8}I_7
    \right).
\end{align}
The control signal sent to the simulation includes percent throttle $u$, steering wheel rate $\dot{\delta}$ in radians per second.
A watermark with covariance 
\begin{align}
    \Sigma_e=
    \begin{bmatrix}
    0.015&0\\0&0.015
    \end{bmatrix}
\end{align}
was added to the control input at each step.
The matrix normalizing factor and the auto-correlation normalizing factor were generated from 200 realizations using \eqref{eq:vn_approx} and  \eqref{eq:gn_approx}.
The window size of 21 steps ($\ell=20$) was used for the statistical tests.
For this window size, a threshold of 181.94 was used based on a false alarm rate of 0.002 for the un-attacked trials.

To generate a replay attack, the measurement signal from one run is recorded and then played back when the simulation is run for a separate realization.
Since an attack need not start at the beginning, we chose to start the attack 50s after the start of the simulation.
Furthermore, since the initial replayed measurement may be inconsistent with what is expected given the current observed state of the system, the attacked measurement instead was linearly interpolated between the true measurement and the replayed measurement over the course of 0.15s.

In practice, an autonomous vehicle would respond to the detection of an attack.
We instead allowed the vehicle to continue normal operation up to a certain distance from the desired trajectory.
This allows us to illustrate the results of a replay attack on an autonomous vehicle.

The results of these simulations can be seen in Figure \ref{fig:carsim}.
The left side of the figure shows the results of the replay attack on our high fidelity car model.
The right side of the figure shows the ability of LTV Dynamic Watermarking to detect these attacks.
Note, despite our attempt to smooth the transition to the replayed attack the negative log likelihood has a spike immediately following the start of the attack at 50 s.
Moreover, the negative log likelihood continues to exceed the threshold as the attack continues and the transient effect of the transition diminishes.

\section{Experimental Results}\label{sec:experimental}
This section further illustrates the effectiveness of LTV Dynamic Watermarking on a 1/10 scale autonomous rover.
For the experiment, the rover completes a lap around a track consisting of several turns and changes in velocity.
The track has a length of 38.8 m and the rover travels at speeds up to 1.8 m/s.
This is accomplished using a LQR and a linearized rover model.
The measurement signal at step $n$ includes the ground plane coordinates $(x_{c,n},y_{c,n})$ in meters, heading $\psi_n$ in radians, angular velocity $\dot{\psi}_n$ in radians per second, and longitudinal velocity $v_{1,n}$ in meters per second.
The ground plane coordinates and heading are measured using a motion capture system, the angular velocity is measured by an imu, and the longitudinal velocity is measured by the motor controller.
The control signal includes a desired speed in meters per second and a steering angle in radians.
A watermark with covariance
\begin{align}
    \Sigma_e=
    \begin{bmatrix}
    0.02&0\\0&0.005
    \end{bmatrix}
\end{align}
was added to the control input at each step.
The matrix normalizing factor and the auto-correlation normalizing factor were generated from 100 experimental runs using \eqref{eq:vn_approx} and \eqref{eq:gn_approx}. 
The window size of 15 steps ($\ell=14$) was used for the statistical tests.
For this window size, a threshold of 175.28 was used based on a false alarm rate of 0.002 for the un-attacked trials.

Implementation of the replay attack was done in the same fashion as was done in simulation except the attack was initiated at 15 s.
For safety purposes, the rover is remotely stopped when the attack causes it to leave the track area.

The results of these experiments can be seen in Figure \ref{fig:rover}.
Similar to the simulated results, the left side of the figure shows the results of the replay attack on the 1/10 scale autonomous rover.
Furthermore, the right of the figure shows the ability of LTV Dynamic Watermarking to detect these attacks.
Note, the transition to the replayed measurements has a lesser effect on the negative log likelihood.
Nonetheless, the negative log likelihood continues to exceed the threshold as the attack continues ensuring detection.
\section{Conclusion}\label{sec:conclusion}
This paper relaxes the requirement for inputs to be visible in a single step and derives an auto-correlation normalizing factor for LTV Dynamic Watermarking.
The new normalizing factor is proven to remove the auto-correlation of the residuals to solidify the statistical background of the implementable test.
The effectiveness of the new normalizing factor is shown using an example.
Furthermore, this paper provides proof of concept for LTV Dynamic Watermarking using both a high-fidelity car model in \carsim and a 1/10 scale autonomous rover.
In each case, a replay attack is implemented in the middle of a trajectory following task, and LTV Dynamic Watermarking is shown to quickly detect the attack in a repeatable fashion.


\renewcommand{\bibfont}{\normalfont\small}
{\renewcommand{\markboth}[2]{}
\printbibliography}

@book{Kollo2005,
    title = {{Advanced Multivariate Statistics with Matrices}},
    year = {2005},
    author = {Kollo, Tonu and Rosen, Dietrich von},
    pages = {237--238},
    publisher = {Springer},
    isbn = {9781402034183}
}

@inproceedings{Humphreys2008,
    title = {{Assessing the Spoofing Threat: Development of a Portable GPS Civilian Spoofer}},
    year = {2008},
    booktitle = {Proceedings of the 21st International Technical Meeting of the Satellite Division of The Institute of Navigation},
    author = {Humphreys, Todd E. and Ledvina, Brent M and Tech, Virginia and Psiaki, Mark L and Hanlon, Brady W O and Kintner, Paul M},
    number = {September 2008},
    pages = {2314--2325},
    isbn = {9781605606897}
}

@inproceedings{Hespanhol2017,
    title = {{Dynamic Watermarking for General LTI Systems}},
    year = {2017},
    booktitle = {56th IEEE Conference on Decision and Control (CDC)},
    author = {Hespanhol, Pedro and Porter, Matthew and Vasudevan, Ram and Aswani, Anil},
    isbn = {9781509028726},
    arxivId = {1703.07760}
}

@article{Satchidanandan2017,
    title = {{Dynamic Watermarking: Active Defense of Networked Cyber-Physical Systems}},
    year = {2017},
    journal = {Proceedings of the IEEE},
    author = {Satchidanandan, Bharadwaj and Kumar, P. R.},
    number = {2},
    pages = {219--240},
    volume = {105},
    isbn = {00189219 (ISSN)},
    doi = {10.1109/JPROC.2016.2575064},
    issn = {00189219},
    arxivId = {1606.08741}
}

@article{Greenberg_2015,
    title = {{Hackers Cut a Corvette’s Brakes Via a Common Car Gadget}},
    year = {2015},
    journal = {Wired.com},
    author = {Greenberg, Andy},
    pages = {},
    url = {http://www.wired.com/2015/08/hackers-cut-corvettes-brakes-via-common-car-gadget/}
}

@article{Greenberg2015,
    title = {{Hackers Remotely Kill a Jeep on the Highway—With Me in It | WIRED}},
    year = {2015},
    journal = {Wired.com},
    author = {Greenberg, Andy},
    pages = {},
    url = {http://www.wired.com/2015/07/hackers-remotely-kill-jeep-highway/}
}

@article{Abrams2008,
    title = {{Malicious Control System Cyber Security Attack Case Study - Maroochy Water Services, Australia}},
    year = {2008},
    journal = {MITRE},
    author = {Abrams, Marshall and Weiss, Joe}
}

@inproceedings{Dominic2016,
    title = {{Risk Assessment for Cooperative Automated Driving}},
    year = {2016},
    booktitle = {Proceedings of the 2nd ACM Workshop on Cyber-Physical Systems Security and Privacy - CPS-SPC '16},
    author = {Dominic, Derrick and Chhawri, Sumeet and Eustice, Ryan M. and Ma, Di and Weimerskirch, André},
    pages = {47--58},
    isbn = {9781450345682},
    doi = {10.1145/2994487.2994499}
}

@inproceedings{Mo2009,
    title = {{Secure Control Against Replay Attacks}},
    year = {2009},
    booktitle = {47th Annual Allerton Conference on Communication, Control, and Computing},
    author = {Mo, Yilin and Sinopoli, Bruno},
    pages = {911--918},
    isbn = {9781424458714},
    doi = {10.1109/ALLERTON.2009.5394956}
}

@incollection{Gerla2015,
    title = {{Securing the Future Autonomous Vehicle: A Cyber-Physical Systems Approach}},
    year = {2015},
    booktitle = {Securing Cyber-Physical Systems},
    author = {Gerla, Mario and Reiher, Peter},
    chapter = {7},
    pages = {197--220},
    publisher = {CRC Press}
}

@article{Amoozadeh2015,
    title = {{Security Vulnerabilities of Connected Vehicle Streams and Their Impact on Cooperative Driving}},
    year = {2015},
    journal = {IEEE Communications Magazine},
    author = {Amoozadeh, Mani and Raghuramu, Arun and Chuah, Chen Nee and Ghosal, Dipak and Michael Zhang, H. and Rowe, Jeff and Levitt, Karl},
    number = {6},
    pages = {126--132},
    volume = {53},
    isbn = {0163-6804 VO - 53},
    doi = {10.1109/MCOM.2015.7120028},
    issn = {01636804},
    pmid = {19221721}
}

@inproceedings{Hespanhol2018,
    title = {{Statistical Watermarking for Networked Control Systems}},
    year = {2018},
    booktitle = {2018 Annual American Control Conference (ACC)},
    author = {Hespanhol, Pedro and Porter, Matthew and Vasudevan, Ram and Aswani, Anil},
    arxivId = {1709.08617},
    pages = {5467--5472}
}

@article{Langner2011,
    title = {{Stuxnet: Dissecting a Cyberwarfare Weapon}},
    year = {2011},
    journal = {IEEE Security and Privacy},
    author = {Langner, Ralph},
    number = {3},
    pages = {49--51},
    volume = {9},
    isbn = {1540-7993},
    doi = {10.1109/MSP.2011.67},
    issn = {15407993}
}

@article{Lee2016AnalysisGrid,
    title = {{Analysis of the Cyber Attack on the Ukrainian Power Grid}},
    year = {2016},
    journal = {Electricity Information Sharing and Analysis Center (E-ISAC)},
    author = {Lee, Rober M. and Assante, Michael J. and Conway, Tim}
}

@book{anderson2003,
  title={An Introduction to Multivariate Statistical Analysis},
  author={Anderson, T.W.},
  isbn={9780471360919},
  lccn={20234317},
  series={Wiley Series in Probability and Statistics},
  year={2003},
  publisher={Wiley}
}

@book{grimmett2001probability,
  title={Probability and Random Processes},
  author={Grimmett, Geoffrey and Stirzaker, David},
  year={2001},
  edition = {3},
  publisher={Oxford university press},
  isbn={0198572220}
}

@article{Porter2019LTV,
    title = {{Detecting Generalized Replay Attacks via Time-Varying Dynamic Watermarking}},
    year = {2019},
    Author = {Porter, Matthew and Hespanhol, Pedro and Aswani, Anil and Johnson-Roberson, Matthew and Vasudevan, Ram},
    journal={arXiv preprint arXiv:1909.08111}
}

@article{hespanhol2019sensor,
  title={{Sensor Switching Control Under Attacks Detectable by Finite Sample Dynamic Watermarking Tests}},
  author={Hespanhol, Pedro and Porter, Matthew and Vasudevan, Ram and Aswani, Anil},
  journal={arXiv preprint arXiv:1909.00014},
  year={2019}
}

@inproceedings{Porter2019,
    author={M. Porter and A. Joshi and P. Hespanhol and R. Vasudevan and A. Aswani}, 
    booktitle={2019 Annual American Control Conference (ACC)},
    title={{Simulation and Real-World Evaluation of Attack Detection Schemes}}, 
    year={2019}, 
    volume={}, 
    number={}, 
    pages={551--558},
    doi={}, 
    ISSN={}, 
}

@inproceedings{Cao2019,
  title={{Adversarial Sensor Attack on LiDAR-based Perception in Autonomous Driving}},
  author={Yulong Cao and Chaowei Xiao and Benjamin Cyr and Yimeng Zhou and Won Park and Sara Rampazzi and Qi Alfred Chen and Kevin Fu and Zhuoqing Morley Mao},
  booktitle={Proceedings of the 26th ACM Conference on Computer and Communications Security (CCS'19)},
  year={2019},
  address = {London, UK}
}

@InProceedings{Eykholt2018,
author = {Eykholt, Kevin and Evtimov, Ivan and Fernandes, Earlence and Li, Bo and Rahmati, Amir and Xiao, Chaowei and Prakash, Atul and Kohno, Tadayoshi and Song, Dawn},
title = {{Robust Physical-World Attacks on Deep Learning Visual Classification}},
booktitle = {The IEEE Conference on Computer Vision and Pattern Recognition (CVPR)},
year = {2018}
}

@article{petit2015remote,
  title={{Remote Attacks on Automated Vehicles Sensors: Experiments on Camera and Lidar}},
  author={Petit, Jonathan and Stottelaar, Bas and Feiri, Michael and Kargl, Frank},
  journal={Black Hat Europe},
  year={2015}
}

@INPROCEEDINGS{Weerakkody2017exposing,
author={S. {Weerakkody} and O. {Ozel} and P. {Griffioen} and B. {Sinopoli}},
booktitle={2017 IEEE Conference on Control Technology and Applications (CCTA)},
title={{Active Detection for Exposing Intelligent Attacks in Control Systems}},
year={2017},
volume={},
number={},
pages={1306-1312},
doi={10.1109/CCTA.2017.8062639},
ISSN={}
}

@INPROCEEDINGS{Weerakkody2016info,
author={S. {Weerakkody} and B. {Sinopoli} and S. {Kar} and A. {Datta}},
booktitle={55th IEEE Conference on Decision and Control (CDC)},
title={{Information Flow for Security in Control Systems}},
year={2016},
volume={},
number={},
pages={5065-5072},
doi={10.1109/CDC.2016.7799044},
ISSN={}
}

@ARTICLE{Shuai2014,
author={Z. {Shuai} and H. {Zhang} and J. {Wang} and J. {Li} and M. {Ouyang}},
journal={IEEE Transactions on Vehicular Technology},
title={{Combined AFS and DYC Control of Four-Wheel-Independent-Drive Electric Vehicles over CAN Network with Time-Varying Delays}},
year={2014},
volume={63},
number={2},
pages={591-602},
doi={10.1109/TVT.2013.2279843},
ISSN={},
}

\appendix 
\begin{proof} (Theorem \ref{thm:alpha})
Assume that $\alpha$ is equal to 0.
Then \eqref{eq:new_ltv_watermark_correlation_test} holds by the same reasoning as for the proof of the original theorem \cite[Theorem III.7]{Porter2019LTV}.

Now assume that \eqref{eq:ltv_watermark_correlation_test} holds. 
Rearranging (\ref{eq:ltv_watermark_correlation_test}) using \eqref{eq:output_equation}, \eqref{eq:observer_error_update_full}, and \eqref{eq:attack_definition} results in 
\begin{align}
    \plim{i}~&\spaceSave\frac{1}{i}\sum_{n=0}^{i-1} V_n(C_n\hat{x}_n-y_n)e_{n-\kappa}^\T=\nonumber\\
    &=\plim{i}~\spaceSave\frac{1}{i}\sum_{n=0}^{i-1} V_n(C_n\delta_n-(1+\alpha)z_n+\nonumber\\
    &\hspace{0.8in}-\alpha C_nx_n-C_n\xi_n-\zeta_n)e_{n-\kappa}^\T.\label{eq:alpha_C1_first_expand}
\end{align}
Note,
\begin{align}
    \plim{i}~\spaceSave\frac{1}{i}\sum_{n=0}^{i-1} V_n(-(1+\alpha)z_n-C_n\xi_n-\zeta_n)e_{n-\kappa}^\T=0_{r\times q}\label{eq:alpha_measure_noise_cancel}
\end{align}
by \citet[Corrolary A.7.]{Porter2019LTV} since $z_n$, $\zeta_n$, $\xi_n$ and $e_{n-\kappa}$ are mutually independent and satisfy the necessary auto-correlation bound.
Then by \citet[Theorem A.4.]{Porter2019LTV} we can cancel these terms resulting in 
\begin{align}
    \plim{i}~&\spaceSave\frac{1}{i}\sum_{n=0}^{i-1} V_n(C_n\hat{x}_n-y_n)e_{n-\kappa}^\T=\nonumber\\
    &=\plim{i}~\spaceSave\frac{1}{i}\sum_{n=0}^{i-1} V_n(C_n\delta_n-\alpha C_nx_n)e_{n-\kappa}^T.\label{eq:alpha_C1_first_cancel}
\end{align}
Expanding $x_n,\delta_n$ in \eqref{eq:alpha_C1_first_cancel} by $\kappa+1$ steps using \eqref{eq:state_update} and \eqref{eq:observer_error_update_full} then collecting all terms that do not depend on $e_{n-\kappa-1}$ and denoting them $a_n$ results in 
\begin{align}
    &\plim{i}~\spaceSave\frac{1}{i}\sum_{n=0}^{i-1} V_n(C_n\hat{x}_n-y_n)e_{n-\kappa}^\T=\nonumber\\
    &~=\plim{i}~\spaceSave\frac{1}{i}\sum_{n=0}^{i-1} V_n\bigg(a_n+\nonumber\\
    &~-\alpha\sum_{j=0}^{\kappa-1} M_{j,n}C_{n-j}\AK_{(n-j-1,n-\kappa+1)}B_{n-\kappa}e_{n-\kappa}\bigg) e_{n-\kappa}^\T.\label{eq:alpha_C1_second_expand}
\end{align}
where $M_{j,n}\in\mathbb{R}^{r\times r}$ is a bounded linear transform due to the dynamics being bounded and $\kappa$ being finite.
Moreover, $M_{0,n}=I_r$ and due to our choice of $\kappa$ terms for $j>0$ can be cancelled by \citet[Theorem A.4.]{Porter2019LTV} since they converge to $0_{q,r}$ by \citet[Corollary a.7.]{Porter2019LTV} resulting in 
\begin{align}
    &\plim{i}~\spaceSave\frac{1}{i}\sum_{n=0}^{i-1} V_n(C_n\hat{x}_n-y_n)e_{n-\kappa}^\T=\nonumber\\
    &~=\plim{i}~\spaceSave\frac{1}{i}\sum_{n=0}^{i-1} V_n\bigg(a_n+\nonumber\\
    &~-\alpha C_{n}\AK_{(n-1,n-\kappa+1)}B_{n-\kappa}e_{n-\kappa}\bigg) e_{n-\kappa}^\T.\label{eq:alpha_C1_third_expand}
\end{align}
Then by \citet[Corollary A.7.]{Porter2019LTV} we have that 
\begin{align}
    \plim{i}~&\spaceSave\frac{1}{i}\sum_{n=0}^{i-1}-\alpha V_nC_n\AK_{(n-1,n-\kappa+1)}B_{n-\kappa}\times\nonumber\\
    &\quad\times(e_{n-\kappa}e_{n-\kappa}^\T-\Sigma_e)=0_{q\times r}.
\end{align}
Therefore by \citet[Theorem A.4.]{Porter2019LTV} we have
\begin{align}
    &\plim{i}~\spaceSave\frac{1}{i}\sum_{n=0}^{i-1} V_n(C_n\hat{x}_n-y_n)e_{n-\kappa}^\T=\nonumber\\
    &~=\plim{i}~\spaceSave\frac{1}{i}\sum_{n=0}^{i-1} V_n a_n e_{n-\kappa}^\T+\nonumber\\
    &\qquad-\alpha V_nC_n\AK_{(n-1,n-\kappa+1)}B_{n-\kappa}\Sigma_e.\label{eq:alpha_C1_second_cancel}
\end{align}
Note, that all elements of
\begin{align}
    V_n a_n e_{n-\kappa}^\T
\end{align}
are distributed symmetrically about 0 for all $n\in\mathbb{N}$ since $a_n$ is a zero mean Gaussian random vector.
Consider an element of \eqref{eq:alpha_C1_second_cancel} for which the corresponding element in
\begin{align}\label{eq:C1_non_converge}
    \spaceSave\frac{1}{i}\sum_{n=0}^{i-1} V_nC_n\AK_{(n-1,n-\kappa+1)}B_{n-\kappa}\Sigma_e
\end{align}
does not converge.
For each $i$, the probability that the matrix element in \eqref{eq:alpha_C1_second_cancel} is farther away from 0 than the corresponding element in \eqref{eq:C1_non_converge} is at least $0.5$.
Therefore the element cannot converge in probability to 0 completing the proof.
\end{proof}

\end{document}